\documentclass[11pt]{article}
\usepackage{amsthm,amsfonts,amsmath,amssymb}
\include{PDF}
\begin{document}

\title{\textbf{The Quantitative Characterization of Finite Simple Groups}
 \thanks{The project is supported in part by the
Natural Science Foundation of China (No. 11171364).
\newline
       \hspace*{0.5cm} \scriptsize\emph{E-mail address:}
       wujieshi43@yahoo.com.cn}
\author{ Wujie Shi\\
Department of Mathematics and Statistics, Chongqing
University\\ of Arts and Sciences, Chongqing 402160, China\\and\\
School of Mathematics, Suzhou  University , Suzhou\\ Jiangsu
215006, China}\\
\textbf{\normalsize In celebration of the 80 birthday of Professor John
Thompson}}

\date{ }
\maketitle

\date{ }
\maketitle
\begin{abstract}
\indent In this report we summarize this work, all finite simple groups $G$ can determined
uniformly using their orders $|G|$ and the set $\pi_e(G)$ of their element orders. \\


 \end{abstract}
\section{Introduction}

Group theory is an important branch of mathematics. It has a wide
range of applications in other mathematics, physics, chemistry, and
other fields. Completed (announced) in 1980 , the classification theorem of
finite simple groups, is one of the most important mathematical
achievements of the 20th century. It is long time to prove this
theorem (D. Gorenstein called it the "Thirty Years War"[1]),
participants in many countries in hundreds of group theory
scientist, many articles (close to 500), length (more than 15,000
pages), is unprecedented in the history of mathematics, with
landmark significance.\\

 For a finite group, the order of group and
the element order are two of the most important basic concepts. Let
$G$ be a finite group and $\pi_e(G)$ be the set of element orders in
$G$. In 1987, the author of this paper posed the following conjecture[2]:\\

{\bf Conjecture.} Let $G$ be a group and $M$ a finite simple group. Then $G\cong M$ if
and only if (a) $\pi_e(G) =  \pi_e(M)$, and (b) $|G| = |M|$.\\

That is, for all finite simple groups we may characterize them using
only their orders and the sets of their element orders (briefly,
"two orders").\\

 After I wrote some letters to Prof. John G. Thompson
and reported the above conjecture. Thompson pointed that, "Good luck
with your conjecture about simple groups. I hope you continue to
work on it", "This would certainly be a nice theorem", in his
reply letters. The warmly encouragement prompt us to finish this
work. \\

Moreover, Thompson posed the following problem and conjecture:\\

For each finite group $G$ and each integer $d \geq 1$, let $G(d) =
\{x\in G; x^d = 1\}$. $G_1$ and $G_2$ are of the same order type if
and only if $|G_1(d)| = |G_2(d)|, d = 1,2,\cdots$.\\

{\bf Thompson's Problem (1987).} Suppose $G_1$ and $G_2$ are groups
of the same order type. Suppose also $G_1$ is solvable. Is it true
that $G_2$ is also necessarily solvable?\\

If $G$ is a finite group, set N($G$) = $\{n \in N; G$ has a
conjugacy class $C$ with $|C| = n \}$.\\

{\bf Thompson's Conjecture (1988).} If $G$ and $M$ are of finite
groups and N($G$) = N($M$), and if in addition, $M$ is a non-Abelian
simple group while the center of $G$ is 1, then $G$ and $M$ are
isomorphic.\\

From 1987 to 2003, the authors of [2-8] proved that this conjecture
is correct for all finite simple groups except $B_n$, $C_n$ and
$D_n$ ($n$ even). In the end of 2009, the authors of [9] proved that
this conjecture is correct for $B_n$, $C_n$ and $D_n$ ($n$ even).
Thus, this conjecture is proved and become a theorem, that is, {\bf all
finite simple groups can determined by their
"two orders"}.\\

    {\bf Question 1.}  Find the application for all finite
simple groups can determined by their "two orders".\\

Let $G$ be a finite group and $B(G)$ be Burnside ring of $G$. We
have the following application:\\

{\bf Corollary 1.} Let $G$ be a finite simple group. Then $B(G)$
determines $G$ up to isomorphism.
{\bf Proof.} See [10, Theorem 5.3.].\\

{\bf Question 2.}  Proving this conjecture, can or not independent
on the classification theorem of finite simple groups? For a small
number of nonabelian simple group, for example, $A_5$, we may do it (see [11]). \\

{\bf Question 3.}  Weaken the condition of "two orders",
characterize all finite simple groups.\\

\end{document}